\documentclass[10pt]{elsarticle}
\usepackage{amssymb}
\usepackage{amsmath}
\usepackage{amsfonts}
\usepackage{amsbsy}
\usepackage{amscd}
\usepackage{amsthm}
\usepackage{times}
\usepackage{psfrag}
\usepackage{graphics}
\usepackage{verbatim}
\usepackage{algorithm}
\usepackage{algorithmicx}
\usepackage{algpseudocode}
\usepackage{listings}
\usepackage{color}
\usepackage{colortbl}
\usepackage{graphicx}
\usepackage{caption}
\usepackage{subcaption}
\usepackage{array,supertabular}
\usepackage{tabularx}
\usepackage{booktabs}
\usepackage{multirow}
\usepackage{ulem}
\usepackage{units}
\usepackage{mathabx}
\usepackage{accents}
\usepackage{algorithmicx}
\usepackage{algorithm}
\usepackage{xfrac}
\usepackage{algpseudocode}
\usepackage{rotating}
\usepackage[capitalize]{cleveref}

\crefname{secinapp}{Section}{Sections}
\Crefname{secinapp}{Section}{Sections}

\newlength{\dhatheight}

\newtheorem{theorem}{Theorem}[section]

\newtheorem{definition}[theorem]{Definition}
\newtheorem{remark}[theorem]{Remark}

\def\XXint#1#2#3{{\setbox0=\hbox{$#1{#2#3}{\int}$ }
\vcenter{\hbox{$#2#3$ }}\kern-.57\wd0}}

\newcounter{myalgorithmctr}

\algnewcommand{\LineComment}[1]{\State \(\triangleright\) #1}

\newcommand{\func}[2]{\text{#1}\left(#2\right)}

\newcommand{\len}[1]{\text{len}\left(#1\right)}

\newcommand{\pluseq}{\mathrel{+}=}

\newcommand{\diveq}{\mathrel{/}=}
\newcommand{\mineq}{\mathrel{-}=}

\newcommand{\pr}[1]{\left(#1\right)}
\newcommand{\br}[1]{\left[#1\right]}
\newcommand{\set}[1]{\left\{#1\right\}}
\newcommand{\halfopen}[1]{\left[#1\right)}

\newcommand{\spanof}[1]{\text{span}\set{#1}}
\newcommand{\shape}[1]{\text{shape}\left(#1\right)}

\newcommand{\N}[2]{N_{#1}^{#2}}
\newcommand{\Nof}[3]{N_{#1}^{#2}\left(#3\right)}

\lstdefinestyle{FromFile}{language=Python,
                          frame=single,
                          numbers=left,
                          numberstyle=\tiny,
                          stepnumber=5,
                          numbersep=7pt,
                          numberfirstline=true,
                          abovecaptionskip=2pt,
                          belowcaptionskip=2pt}

\begin{document}

\begin{frontmatter}
\title{A New Algorithm for the Evaluation of Generalized B-splines}

\author[byu1]{Ian D. Henriksen}

\author[byu1]{Emily J. Evans\corref{cor1}}
\ead{ejevans@mathematics.byu.edu}

\author[byu3]{D. C. Thomas}

\cortext[cor1]{Corresponding author}

\address[byu1]{Department of Mathematics,
  Brigham Young University,
  Provo, Utah 84602, USA}

\address[byu3]{Isogeometrx,
  Mapleton, Utah 84664, USA}

\begin{abstract}
In this paper we present a method for direct evaluation of generalized B-splines (GB-splines) via the local representation of these curves as piecewise functions.  To accomplish this we introduce a local structure that makes GB-spline curves more amenable to the techniques used in constructing bases of higher degree.
This basis is used to perform direct computation of piecewise representation of GB-spline bases and curves.
\end{abstract}


\end{frontmatter}

\section{Introduction}
In computer aided geometric design (CAGD) the use of non-uniform rational B-splines (NURBS) as the basis for design is prevalent. The use of NURBS as the basis for geometric design is not without issues, however.  First, NURBS cannot represent certain transcendental curves, many of which such as the helix and cycloid are used in design.  Second, NURBS requires the use of weights to describe certain curves, the selection of which have no geometric meaning.  Lastly the parametrization of conic sections does not correspond to the natural arc-length parametrization.  Much research has been done in the computer aided design community to develop alternative technologies to the standard NURBS technology.  Generalized B-splines (GB-splines) are one such technology that has received increased attention in recent years.

GB-splines are a generalization of B-splines that resolve some of the fundamental problems with the use of NURBS.
Rather than spanning the spaces of piecewise polynomials spanned by traditional B-spline curves, on each interval $\halfopen{t_i, t_{i+1}}$ in the given knot vector $T$, they span the spaces $\{1, t, \dots, t^{p-2}, u_i^{\br{p-1}}, v_i^{\br{p-1}}\}$ where $u_i^{\br{p-1}}$ and $v_i^{\br{p-1}}$ are $p-1^{th}$ integrals of arbitrary functions forming a Chebyshev space over $[t_i, t_{i+1}]$.
Because of their ability to span more general classes of functions, GB-splines allow exact representation of polynomial curves, helices and conic sections  using control point representations that are intuitive and natural to designers~\cite{shapepreserv}.  GB-splines possess all of the fundamental properties of B-splines and NURBS that are important for design and analysis such as local linear independence, degree-elevation and partition of unity.  In addition to the geometric advantages of using GB-splines over NURBS, GB-spines also behave similarly to B-splines with respect to differentiation and integration. This similarity in behavior is especially beneficial when relevant properties of the continuous problem must be transferred to the discrete problem~\cite{electromag}.

In 1999 Ksasov and Sattayatham \cite{Ksasov1999} demonstrated a variety of the properties of GB-splines.
In 2005 Costantini et al. \cite{costantini2005} studied generalized Bernstein bases of this form in greater detail.
In 2008, Wang et al. \cite{ue_spline_original} introduced unified extended splines (UE-splines), a subclass of GB-splines and demonstrated that this new class of splines contains several other classes of generalized splines.
In 2011, Manni et al. \cite{iga_gb_splines} proposed that GB-splines be used for isogeometric analysis.
In \cite{exponential_subdivision}, Romani successfully applied the techniques from \cite{subdivision_book} to form subdivision methods that allow for the approximation of UE-splines via a limit of control meshes successively refined by a non-stationary subdivision scheme. In~\cite{quasiinterpolation}, quasi-interpolation was performed in isogeometric analysis using GB-splines.  Finally, GB-splines have also been used in the context of T-meshes~\cite{tmeshes}.

From the usual recursive definition of GB-splines the only effective means of evaluating GB-splines is either through recursive numeric integration, or through symbolic computation of indefinite integrals. Recursive numeric integration is very costly for all but the lowest degrees of splines, whereas
symbolic computation, while effective, can be unwieldy for numeric computation.
In order to address these difficulties, we present a more direct method of computing values on GB-spline curves, using local representations.

\subsection{Structure and content of the paper}
\label{sec:content}
In~\cref{sec:GBRev} the GB-splines are reviewed and appropriate notational
conventions are introduced. An algorithm for their direct evaluation is introduced 
in~\cref{sec:aegs}. 
In ~\cref{sec:conclusion} we draw conclusions. 

\section{A review of generalized B-splines}
\label{sec:GBRev}
Generalized B-splines (GB-splines) were introduced in \cite{ue_spline_original}, and span spaces of the form $\set{1, t, \dots, t^{p-2}, u\pr{t}, v\pr{t}}$ where $u$ and $v$ are more general functions defined over each interval in a knot vector. GB-splines retain most of the desirable properties of B-splines and unify a variety of other spline types such as UE-splines, trignometric splines, exponential splines, etc.    The primary advantages of GB-spline curves over traditional B-splines is that they allow for the exact representation of certain geometric curves and surfaces, like circles, hyperbolas, spheres, and hyperboloids,that cannot be well-represented by polynomial splines.
Before before formally defining a GB-spline we present some preliminary definitions.
\begin{definition} \label{knot_vector_definition}
A knot vector, $T$, is a nondecreasing vector of real numbers.  
\end{definition}

\begin{definition}
A set of $\ell$ linearly independent functions are said to form a Chebyshev space over an interval $I$ if any nonzero function in their span has at most $\ell-1$ roots in that interval.
\end{definition}
\begin{definition}
Given a knot vector $T$, and functions $u_i$ and $v_i$ forming a Chebyshev space on each $\br{t_{i}, t_{i+1}}$ of nonzero length such that $u\pr{0} = v\pr{1} = 1$ and $v\pr{0} = u\pr{1} = 0$, we will refer to the sets of functions $u_i$ and $v_i$ as the knot functions over $T$.
\end{definition}
With these definitions we are prepared to define a GB-spline.
\begin{definition} \label{gbspline_definition}
Given a degree $p$ and a knot vector $T$ of length $m$ with corresponding knot functions $u_i$ and $v_i$.
Define the $i^{th}$ GB-spline basis function of degree $p$, denoted by $\N{i}{p}$ as follows:\\

\noindent Define the degree 1 GB-spline basis function as:
\[\Nof{i}{1}{t} = \begin{cases} u_{i}(t) & t \in \halfopen{t_i, t_{i+1}} \\ v_{i+1}(t) & t \in \br{t_{i+1}, t_{i+2}} \\ 0 & \text{otherwise} \end{cases},\]
For $p\geq 1$ define
\[\delta_i^p = \int_{t_i}^{t_{i+p+1}} \Nof{i}{p}{s} ds.\]
For $p\geq 1$ define $\Phi_i^p(t)$ as 
\[\Phi_i^p\pr{t} = \begin{cases} \frac{\int_{t_i}^t \Nof{i}{p}{s} ds}{\delta_i^p} & \text{if $\delta_i^p \neq 0$,}\\
0 & \text{if $ \delta_i^p = 0$ and $t < t_{i+p+1}$,}\\
1 &  \text{if $ \delta_i^p = 0$ and $t \geq t_{i+p+1}$.} \end{cases}\]
For $p > 1$, define
\[\Nof{i}{p}{t} = \Phi_i^{p-1}\pr{t} - \Phi_{i+1}^{p-1}\pr{t}.\]

In addition, if $t_{m-p-1} = \dots = t_{m-1}$ and $t_{m-p-2} \neq t_{m-p-1}$ (that is, if the last $p$ knots are repeated, and the last basis function is nonzero), define $N_{m-p-2}^p\pr{t_{m-1}} = 1$.
\end{definition}

\begin{definition}
Given a degree $p >1$, and a knot vector $T$ of length $m$ with a corresponding set of knot functions, a degree $p>1$, and $m-p-1$ control points $a_i$, define the corresponding GB-spline curve $f(t)$ as
\[f\pr{t} = \sum_{i=0}^{m-p-2} a_i \Nof{i}{p}{t}\]
for $t \in \br{t_{p}, t_{m-p-1}}$.
\end{definition}
\begin{remark}
In these definitions, and in the algorithms presented later in the manuscript we used zero-based indexing on the basis functions.  That is to say instead of indexing the basis functions from $1, \ldots, m$, we instead index the basis functions from $0, \ldots, m-1.$
\end{remark}

%
%
\begin{definition}
Given a knot vector $T$ with corresponding sets of knot functions $u_i$ and $v_i$, define
\[V_{i}^{p} = \spanof{1, \pr{t-t_i}, \dots, \pr{t-t_i}^{p-2}, u^{\br{p-1}}\pr{t}, v^{\br{p-1}}\pr{t}}\]
Where $u^{\br{p-1}}$, and $v^{\br{p-1}}$ are the $\pr{p-1}^{th}$ indefinite integrals of $u$ and $v$ respectively.
\end{definition}

GB-splines have the following important properties:
\begin{itemize}
\item B-splines are GB-splines~\cite{ue_spline_original}.
\item The support of $\N{i}{p}$ is zero outside the interval $\br{t_{i}, t_{i+p+1}}$.
\item Partition of unity.
\item GB-spline basis functions are linearly independent and positive on the interior of their supports~\cite{Ksasov1999, iga_gb_splines}.
\item GB-spline curves are variation diminishing~\cite{Ksasov1999, iga_gb_splines}.
\item A GB-spline over an open knot vector $T$ with no degenerate basis functions in the corresponding spline basis interpolates its endpoints.
\item A GB-spline curve $f(t)$ of degree $p$ over a knot vector $T$ with knot functions $u_i$ and $v_i$ has the following properties:
\begin{itemize}
	\item The GB-spline basis restricted to each interval lies in $V_i^p$. \label{gb_span}
	\item Each GB-spline is $C^{p-k}$ at each of the knots in the knot interval where $k$ is the number of times a knot is repeated. \label{gb_knot_smoothness}
	\item Each GB-spline is at least $C^{p+r-1}$ for each point $t$ not in its knot vector, where $r$ is the minimum continuity of the knot functions over the interval in the knot vector that contains $t$. \label{gb_interior_smoothness}
\end{itemize}
\item Where it exists, the derivative of a GB-spline basis function is given by
\[\pr{\N{i}{p}}'\pr{t} = \frac{\Nof{i}{p-1}{t}}{\delta_{i}^{p-1}} - \frac{\Nof{i+1}{p-1}{t}}{\delta_{i+1}^{p-1}}.\]
\item Where it exists, the derivative of a GB-spline curve $f(t)$ with control points $a_0, \dots, a_{n}$ is given by
\[\sum_{i=0}^{n+1} a_i \frac{\Nof{i}{p-1}{t}}{\delta_{i}^{p-1}} - a_{i+1} \frac{\Nof{i+1}{p-1}{t}}{\delta_{i+1}^{p-1}}\]
with $a_{-1}$ and $a_{n+1}$ defined to be $0$.
\end{itemize}
\section{An algorithm for evaluation of GB-splines}
\label{sec:aegs}
Although it is clear from Definition~\ref{gbspline_definition} why many of the properties of GB-spline curves are true, it does not provide for a simple means of evaluation. From the definition the only effective means of evaluating GB-splines are either recursive numeric integration, or symbolic computation of indefinite integrals. Recursive numeric integration is very costly for all but the lowest degrees of splines, whereas
symbolic computation, while effective, can be unwieldy for numeric computation.
In order to address these difficulties, we present a more direct method of computing values on GB-spline curves.

Given that each basis function lies in the space $V_i^p$, we may introduce a local representation of each basis function in terms of the functions spanning the space $V_i^p$.  Since the recursive integrals must be computed, we would like for these local representations to be more amenable to integration. These requirements introduce several possible choices for the local representations of the splines, but we will use the local representation given by $u_i^{[p-1]}$, $v_i^{[p-1]}$, and an additional polynomial term of degree $p-2$. By virtue of the linear independence of $u_i^{[p-1]}$ and $v_i^{[p-1]}$ from all other polynomial terms this is a valid choice of basis.

\subsection{Local Representations: Knot Functions and Polynomials}
\begin{definition} \label{gbspline_standard_recurrence}
Given a degree $p$ and a knot vector $T$ of length $m$ with corresponding knot functions $u_i$ and $v_i$, and since $\N{i}{p}$ lies in the space $V_i^p$, we can represent $\N{i}{p}$ on the $j^{th}$ interval in $T$ as:
\[\Nof{i}{p}{t} = P_{i,j}^p(t) + a_{i,j}^p u_j^{[p-1]}(t) + b_{i,j}^p v_j^{[p-1]}(t)\]
where $P_{i,j}^p$ is a polynomial term and $a_{i,j}^p$ and $b_{i,j}^p$ are constants.
The recurrence stated in Definition \ref{gbspline_definition} can be written as:
\[a_{i,j}^p = \frac{a_{i,j}^{p-1}}{\delta_i^{p-1}} - \frac{a_{i+1,j}^{p-1}}{\delta_{i+1}^{p-1}},\]
\[b_{i,j}^p = \frac{b_{i,j}^{p-1}}{\delta_i^{p-1}} - \frac{b_{i+1,j}^{p-1}}{\delta_{i+1}^{p-1}},\]
and
\[\begin{aligned}P_{i,j}^p (t) &= \Nof{i}{p}{t_j} \\&+ \frac{1}{\delta_i^{p-1}} \left(\int_{t_j}^t P_{i,j}^{p-1}(s) ds - a_{i,j}^{p-1} u_j^{[p-1]}(t_j) - b_{i,j}^{p-1} v_j^{[p-1]}(t_j)\right) \\
	&- \frac{1}{\delta_{i+1}^{p-1}} \left(\int_{t_j}^t P_{i+1,j}^{p-1}(s) ds - a_{i+1,j}^{p-1} u_j^{[p-1]}(t_j) - b_{i+1,j}^{p-1} v_j^{[p-1]}(t_j) \right).
\end{aligned}\]

With the additional stipulation that if $\N{i}{p-1}$ is identically zero and the interval $[t_{i+p}, t_{i+p+1})$ is empty, then an additional 1 is added to $P_{i,j}$ to account for the modified treatment of basis functions that are identically $0$ in Definition \ref{gbspline_definition}.
As before, we also require that, if the last basis function is discontinuous at the end of the $m-p-1$ knot, that it must take a value of $1$ at that knot.
\end{definition}

\begin{remark}
For ease of notation later in the manuscript we name the coefficients $a_{i,j}^p$ and $b_{i,j}^p$ in Definition~\ref{gbspline_standard_recurrence} as the \emph{general function coefficients}.
\end{remark}
\begin{figure}
\begin{minipage}[b]{0.5\linewidth}
\centering
\includegraphics[width=\textwidth]{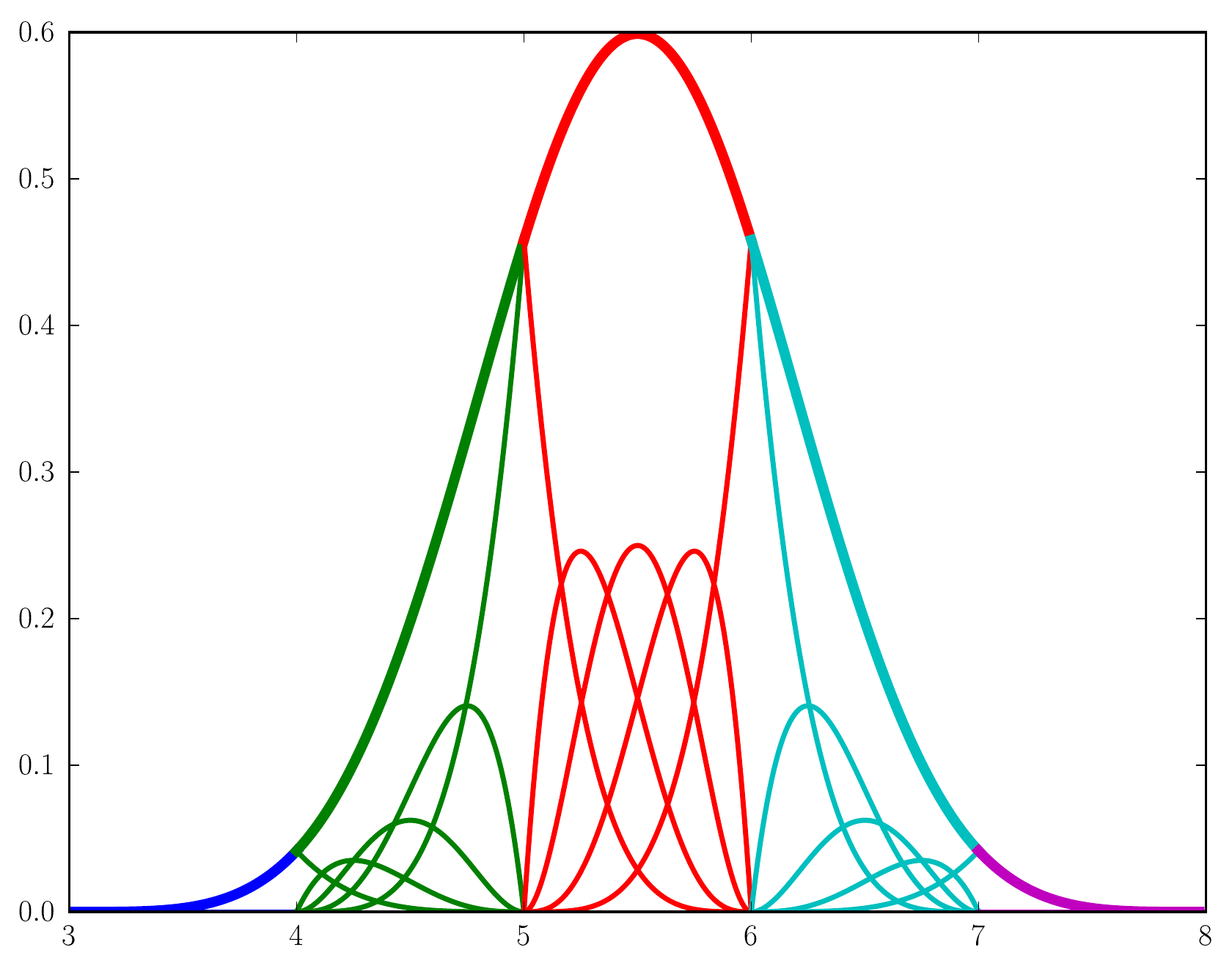}
\end{minipage}
\begin{minipage}[b]{0.5\linewidth}
\centering
\includegraphics[width=\textwidth]{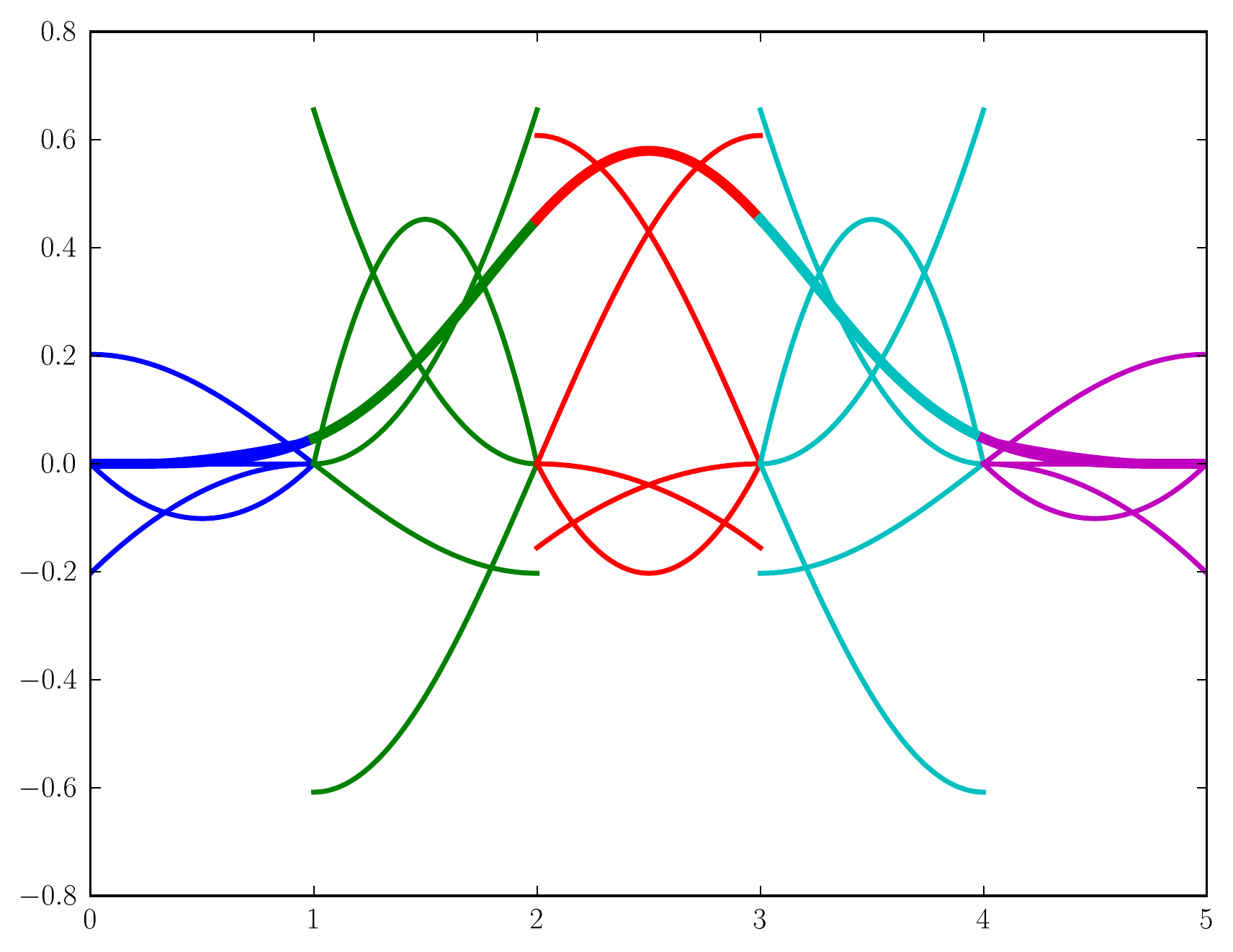}
\end{minipage}
\caption{The local polynomial representation for a uniform B-spline basis function of degree $4$ compared with the local representation for a uniform GB-spline function of degree $4$ defined using trigonometric knot functions.}
\end{figure}

The recurrence relation outlined in Definition \ref{gbspline_standard_recurrence} is not as easy to implement as De Boor's recurrence, however it does make it so that the evaluation of GB-spline curves is no longer tied to symbolic integrals or recursive quadrature.
It makes it clear that the values of $\N{i}{p}$ on the interval $[t_j, t_{j+1})$ depend only on $\Nof{i}{p}{t_j}$, the values of $u_j^{[p-1]}$ and $v_j^{[p-1]}$ at $t$ and $t_j$, and the full set of coefficients for $\N{i}{p-1}$ and $\N{i+1}{p-1}$.
These dependencies can be stated more explicitly.
To evaluate $\N{i}{p}$ at time $t \in [t_j, t_{j+1})$, it is necessary to know the values of the following function values:
\begin{itemize}
\item $u_j^{[p-1]}$ and $v_j^{[p-1]}$ at $t$;
\item the values of all the different $u_{j}, u_j^{[1]}, \dots, u_j^{[p-2]}$ and $v_j, v_j^{[1]}, \dots, v_j^{[p-2]}$ at $t_j$ and $t_j+1$ for each index $j$ corresponding to an interval in the support of $\N{i}{p}$;
\item the values of $u_j^{[p-1]}$ and $v_j^{[p-1]}$ if $j$ is an index corresponding to an interval in the support of $\N{i}{p}$ at $t_j$ if $t_j < t$  and at $t_{j+1}$ when $t_{j+1} < t$.
\end{itemize}

\subsection{An algorithm for evaluation of GB-splines}

To construct a direct algorithm for the evaluation of arbitrary GB-spline curves, we must first determine how best to handle the dependencies between the intervals.
Given that the representation of $\N{i}{p}$ on $[t_{j}, t_{j+1})$ depends on the representations of $\N{i}{p-1}$ and $\N{i+1}{p-1}$ over their supports and the representations of $\N{i}{p}$ over the intervals of its support that lie to the left of $[t_j, t_{j+1})$, it is natural to construct the set of basis functions of each degree using the set of basis functions of the degree one less than the one being computed.
The computation most naturally runs from left to right along each basis function.
Given this structure, the algorithm should operate roughly as follows:
\begin{itemize}

\item Initialize a list of basis functions using the known values for the degree $1$ case.

\item For each degree from $2$ to the desired degree $p$, do the following:
	\begin{itemize}
	
	\item Integrate each polynomial term in the basis.
	
	\item Use the polynomial terms, the general function coefficients, and the values of the indefinite integrals of the knot functions at the points in the knot vector to compute the definite integral of each basis function over its support.
	
	\item Divide the indefinite integrals of the polynomial terms by the definite integral of the basis function they represent.
	
	\item Divide the general function coefficients by the definite integral of the basis function they represent.
	
	\item Compute the differences between the scaled general function coefficients for basis functions whose indices differ by $1$.
	
	\item Compute the differences between the scaled polynomials for basis functions whose indices differ by $1$, adding the constant terms from the general function part to the polynomial.
	
	\item Use the values of these differences to add the value of each basis function over each interval to its polynomial term over each interval.
	
	\item Store these differences between the polynomial and general function terms as the new set of basis functions.
	
	\end{itemize}

\end{itemize}

In practice, the functions that we desire to include in the span of the spline basis may not always satisfy the constraints on the values of the knot functions at the endpoints of each interval.
This can be resolved by taking linear combinations of the original functions on each interval so that the endpoint constraints are satisfied.
This can be taken care of as a part of the algorithm for constructing a basis by taking the needed linear combinations of the integral terms given as input and then, once the local representations of the spline basis have been computed, changing the representations so that they are given with respect to the original functions rather than the computed linear combinations.
In order to ensure the desired properties of a spline basis, the functions used to create the knot functions must still form a Chebyshev space over each corresponding nondegenerate interval in the knot vector.
The matrix
\[\begin{bmatrix}
u_i\pr{t_i} & v_i\pr{t_i} \\
u_i\pr{t_{i+1}} & v_i\pr{t_{i+1}}
\end{bmatrix}\]
must also be invertible (and sufficiently well-conditioned) so that the needed linear combinations can be computed.

In addition, the spline basis constructed will span $u^{[p-1]}$ and $v^{[p-1]}$, not $u$ and $v$.
It is often desirable to construct a basis that spans $C^{p-1}$ functions $\tilde{u}$ and $\tilde{v}$ instead.
To handle that properly, we need only begin the iteration with the knot functions $\tilde{u}^{(p-1)}$ and $\tilde{v}^{(p-1)}$, noting that, after the corresponding numbers of integrals have all been taken, the spline basis will span the desired functions.
To use this approach, it is necessary that there exist linear combinations of $\tilde{u}^{(p-1)}$ and $\tilde{v}^{(p-1)}$ that satisfy the constraints that would normally be imposed on $u$ and $v$.

The recurrence outlined in Definition \ref{gbspline_standard_recurrence} is restated as an algorithm in Algorithm \ref{gb_basis_construction}.
For convenience, each basis is stored as two arrays, the first containing the polynomial terms corresponding to each interval within the support of each basis function and the second containing the corresponding abstract function terms, that is the terms for $u$ and $v$.
Given that the support of each basis function is known, we only include the representation of each basis function on an interval where it will be nonzero.
This shifts the indices for the representation of each basis function, but the structure of the iteration is essentially the same.  The algorithms here will be presented in a form that is independent of the polynomial basis used.

For practical use it is also helpful to follow the convention that the knot functions and polynomials corresponding to each interval are defined on the interval $\br{0, t_{i+1} - t_i}$ and that $t - t_i$ is used as an argument rather than $t$ itself.
This makes it so that, for any given polynomial representation requiring boundaries of definition (Bernstein polynomials, Chebyshev polynomials, etc.), only the lengths of each interval must be passed to the polynomial integration and evaluation routines.

An important consequence of using the piecewise representations for these basis functions is that, once a piecewise representation for a spline curve is created, the only remaining costs of evaluating the function at any given point come from identifying which interval in the knot vector contains the given parameter value, evaluating a polynomial term, and evaluating the terms $u_i^{\br{p-1}}$ and $v_i^{\br{p-1}}$.
No further recursion or integration is needed.

The algorithms will be presented in vectorized form with a particular emphasis on clarity.
A variety of other small optimizations could be added to further remove redundant computations; however the presentation here is meant primarily to provide a clear explanation of the algorithm.
It presents a relatively efficient version of the algorithm, but, for simplicity, redundant computations have not been completely removed.

For clarity within the algorithm and its helper routines, we will now introduce the many variables used throughout this algorithm and its corresponding helper routines.
Throughout the code for this algorithm, the following variables will be used:
\begin{itemize}
\item
$T$ is a $1$-dimensional array containing the knot vector.

\item
\textit{Tlens} is a $1$-dimensional array containing the lengths of the intervals between the knots values in $T$.

\item
\textit{Tvals} an array shaped like \textit{Tlens} containing the lengths of each interval.
 \textit{Tvals} is indexed first by basis function, then by interval within the support of a given basis function.

\item
$p$ is the degree of the desired basis, or of the spline to be evaluated.

\item $n$ is the number of basis functions in a given basis.

\item
\textit{ints} is a 4-dimensional array of shape $\pr{p, \len{t}-1, 2, 2}$ containing the values of the indefinite integrals of the knot functions at the endpoints of each interval.  
The $p^{th}$ integrals are indexed in ascending order along the first axis.
The different intervals within the knot vector areindexed along the second axis.
The different endpoints of each interval are indexed along the third axis.
The different knot functions ($u$ and $v$) are indexed along the last axis with $u$ first.

\item \textit{wints} is a 3-dimensional view into \textit{ints} corresponding to the integral values of a given degree, indexed first by basis function, then by interval within the support of each basis function, then by endpoint, then by the different knot functions.

\item
\textit{polys} is an array containing the coefficients for the polynomial parts of the basis functions in a given basis.
Basis functions are indexed along the first axis and intervals within the support of each basis function are indexed along the second axis.
Here we will assume that the polynomial term over the $i$'th interval is stored in the form $p\pr{t - t_i}$, that is, that the polynomial terms are translated so that the first value taken by the polynomial in each interval is the value of the polynomial at $0$.
This algorithm does not depend on the representation used to store the polynomial terms, but in most cases an array of shape $\pr{n, p+1, p-1}$ containing only the necessary coefficients should suffice.

\item
\textit{pints} is an array containing the integrals of all the polynomial terms in a given array \textit{polys}.
All the axes are indexed the same as the axes in \textit{polys}.
The shape will be the same except that the last axis will be one index longer than the last axis of \textit{polys}.

\item
\textit{genfunc} will be an array containing the coefficients for the general function terms of the basis functions in a given basis.
Basis functions should be indexed along the first axis, intervals within the support of each basis function along the second axis, and the two general function coefficients along the third (with $u$ first, then $v$).
This array will have a shape of $\pr{n, p+1, 2}$.

\item
\textit{scal} will be an array containing the necessary scaling matrices needed to scale \textit{ints} to represent the scaled versions of $u$ and $v$ that satisfy the value constraints at their endpoints and also needed to scale the coefficients in \textit{genfunc} so that they represent the basis functions in terms of the original $u$ and $v$.

\item
\textit{pos} will be an array of boolean values.
The $i$'th entry of \textit{pos} will be true if $\delta_i^{d-1} \neq 0$.

\item
\textit{deltas} will be an array containing the indefinite integrals of each basis function over its support.

\item
\textit{consts} will be an array containing the constant terms to be added to the polynomial terms on each interval.
It will be indexed first by basis function, then by interval within the support of each basis function.
In the recurrence in Definition \ref{gbspline_standard_recurrence}, these are the terms
\[- a_{i,j}^{d-2} u_j^{[d-2]}(t_j) - b_{i,j}^{d-2} v_j^{[d-2]}(t_j)\]

\item
\textit{vals} will be a temporary array used to store outputs of various functions.

\item
$d$ will be a looping variable used in the loop that constructs the basis of each degree from the basis of previous degree.
Here $d$ will be the degree of the basis being constructed.
\textit{dmin} will be equal to $d-1$.

\end{itemize}

Algorithm \ref{gb_basis_construction} shows the primary routine used to compute the local representations of a given basis function.
It contains calls to a variety of auxiliary routines, all of which will be explained here.
Here we include the primary algorithm first so that the reader may understand the general flow of the algorithm and the proper place for each auxiliary routine before handling the many details that are taken care of in the auxiliary routines.
\begin{algorithm}
\caption{Computing the local coefficients for a GB-spline basis}
\begin{algorithmic}[1]
\Procedure{BasisCoefs}{$T, ints, tol=10^{-8}$}
	\LineComment{Initialize $p$, \textit{Tlens}, and \textit{Tvals} and coefficient arrays for a basis of degree $1$.}
	\State $p = \shape{ints}[0]$
	\State $Tlens = T\br{1:} - T\br{:-1}$
	\State $Tvals = \func{Wrap}{Tlens, 2}$
	\State $polys, genfunc = \func{MakeDegreeOne}{\shape{T}[0] - 2}$
	\LineComment{Take linear combinations of the functions with integrals in \textit{ints} so that}
	\LineComment{the resulting linear combinations satisfy constraints on knot functions.}
	\State $ints, scal = \func{ScaleKnotFuncs}{ints, Tlens, tol}$
	\LineComment{Construct each successive set of local coefficients.}
	\For{$d=2$, $d \leq p$}
		\LineComment{Compute the indefinite integrals of all polynomial terms}
		\State $pints = \func{PolyInt}{polys, Tvals}$
		\LineComment{Construct \textit{wints} by wrapping the first axis of \textit{ints} into two new axes.}
		\State $wints = \func{Wrap}{ints\br{d-1}, d}$
		\LineComment{Integrate the current set of basis functions over their supports.}
		\State $deltas, consts = \func{IntegrateSupports}{Tvals, pints, genfunc, wints}$
		\LineComment{Add constant terms from the general function integrals to the \textit{pints}.}
		\State $\func{OffsetConstants}{pints, consts}$
		\LineComment{Compute the indices of the basis functions that are identically $0$.}
		\State $pos = \pr{T\br{d:} - T\br{:-d}} > tol$
		\LineComment{Take the deltas corresponding to positive basis functions.}
		\LineComment{Also reshape the deltas for broadcasting with \textit{pints} and \textit{genfunc}.}
		\State $deltas = deltas\br{pos, None, None}$
		\LineComment{Divide the terms in \textit{pints} and \textit{genfunc} by their corresponding entries in \textit{deltas}.}
		\State $pints\br{pos} \diveq deltas$
		\State $genfunc\br{pos} \diveq deltas$
		\LineComment{Take the differences between neighboring terms in \textit{pints} and \textit{genfunc}.}
		\State $polys, genfunc = \func{OffsetDifferences}{pints, genfunc}$
		\LineComment{Add ones where needed to account for the integral terms of $0$-valued}
		\LineComment{basis functions after the knot value where their support would end.}
		\State $\func{AddOnes}{polys, pos}$
		\LineComment{Compute the set of Tvals for the next basis.}
		\State $Tvals = \func{Wrap}{Tlens, d+1}$
		\LineComment{Add in the constant terms that come from evaluating each basis function}
		\LineComment{at the end of each interval within the knot vector.}
		\State $\func{ConnectBoundaries}{polys, genfunc, wints, Tvals}$
	\EndFor
	\LineComment{Scale the coefficients in \textit{genfunc} so that the basis functions are represented}
	\LineComment{in terms of the general function terms originally represented in \textit{ints}.}
	\State $\func{ScaleGenFuncCoefs}{\func{Wrap}{scal, p+1}, genfunc}$
	\State \Return $polys, genfunc$
\EndProcedure
\end{algorithmic}
\label{gb_basis_construction}
\end{algorithm}

Algorithm \ref{gb_basis_construction} uses the following auxiliary routines:
\begin{itemize}

\item
\textit{Wrap}: A utility function used to convert between indexing by interval to indexing by basis function, then by interval within the support of each basis function.

\item
\textit{MatMul}: A utility function that performs matrix multiplication.

\item
\textit{MakeDegreeOne}: A function that initializes the coefficient arrays for a basis of degree $1$.

\item
\textit{ScaleKnotFuncs}: A function that computes the scaled \textit{ints} and the corresponding array $invs$ containing the scalings.
This function is what changes all the basis functions to be represented in terms of the linear combinations of the functions used to create $invs$ that satisfy the required constraints to be knot functions.

\item
\textit{ScaleGenFuncCoefs}: A function that modifies \textit{genfunc} in-place to change the coefficients to represent the basis functions in terms of the functions used to create $invs$.

\item
\textit{PolyInt}: A function that, given an array of polynomial coefficients with the coefficients indexed along the last axis and another array containing the lengths of the intervals corresponding to each polynomial term, computes the indefinite integrals of all the polynomials.

\item
\textit{PolyVal}: A function that, given an array of polynomial coefficients with the coefficients indexed along the last axis, and an array containing the lengths of the intervals corresponding to each polynomial term, evaluates each polynomial at the corresponding term in an array \textit{vals}.
This function is used only within other auxiliary routines.

\item
\textit{IntegrateSupports}: A function that computes the integral of each basis function over its support

\item
\textit{GenFuncInts}: A helper function called within \textit{IntegrateSupports}.
It computes the portion of the integral of each basis function that comes from the knot functions over each interval in its support.

\item
\textit{OffsetConstants}: A function that modifies \textit{pints} in-place to add in the constant terms that come from the general function integral.
In the polynomial part of the recurrence from Definition \ref{gbspline_standard_recurrence}, this accounts for the terms
\[- a_{i,j}^{p-1} u_j^{[p-1]}(t_j) - b_{i,j}^{p-1} v_j^{[p-1]}(t_j)\]
and
\[a_{i+1,j}^{p-1} u_j^{[p-1]}(t_j) - b_{i+1,j}^{p-1} v_j^{[p-1]}(t_j)\]

\item
\textit{OffsetDifferences}: A function that takes the \textit{pints} and \textit{genfunc} (after each term has been divided by the corresponding $\delta_i$ terms) that correspond to a given basis and computes the differences between consecutive terms.
This returns the differences between both the polynomial and general function terms.
These terms account for all terms in the recurrence from Definition \ref{gbspline_standard_recurrence} with the exception of $\Nof{i}{p}{t_j}$.

\item
\textit{AddOnes}: A function that adds the one terms to \textit{polys} that come from the handling of the integral terms from basis functions that are identically $0$ as defined in Definition \ref{gbspline_definition}.

\item
\textit{ConnectBoundaries}: A function that computes the term $\Nof{i}{p}{t_j}$ for each interval where it is needed and adds it in place to \textit{polys}.

\end{itemize}

\subsection{ Auxiliary Routines }

In this section we discuss the implementations of the auxiliary routines in greater detail.

The helper routine \textit{Wrap} is used to expand a given axis into two axes where each index of the first of the two new axes provides a moving window of the given width along the original axis that is being expanded.
Within this algorithm, this function is used to convert between data that is indexed by interval within the knot vector to data that is indexed by basis function, then by interval within the support of each basis function.

The implementation for \textit{MakeDegreeOne} should also be relatively simple.
Recall that each degree $1$ basis function has the form
\[\begin{cases}
v_i(t) & x \in \halfopen{t_i, t_{i+1}} \\
u_i(t) & t \in \halfopen{t_{i+1}, t_{i+2}}
\end{cases}\]
This means that this function should allocate \textit{polys} as an empty array of shape $\pr{n, 2, 0}$ and allocate \textit{genfunc} as an array of $0$'s of shape $\pr{n, 2, 2}$.
Then it should fill \textit{genfunc} with values such that $genfunc\br{i,j,k}$ is equal to $0$ when $j = k$ and $1$ otherwise.
It should then return \textit{polys} and \textit{genfunc}.

\textit{ScaleKnotFuncs} allows derivatives from functions that do not necessarily satisfy the constraints $u_i\pr{t_i} = v_i\pr{t_{i+1}} = 1$ and $u_i\pr{t_{i+1}} = v_i\pr{t_i} = 0$ for the integral values stored in \textit{ints}.
To do this, it must require that, for each nonempty interval $\halfopen{t_i, t_{i+1}}$, the matrix
\[A_i = \begin{bmatrix}
u_i\pr{t_i} & v_i\pr{t_i} \\
u_i\pr{t_{i+1}} & v_i\pr{t_{i+1}}
\end{bmatrix}\]
be invertible and reasonably well-conditioned.
It is still required that $u_i$ and $v_i$ form a Chebyshev space.

Since the algorithm for constructing the basis coefficients relies on each $u_i$ and $v_i$ satisfying the constraints on its values at the endpoints of each interval in the knot vector, we must compute the linear combinations of $u$ and $v$ that satisfy the value constraints at each endpoint.
Since matrix multiplication can be seen as using the columns of the matrix on the right as coefficients for linear combinations of the columns of the matrix on the left, we see that the matrix $B_i$ with the desired coefficients for the linear combinations must satisfy the equation $A_i B_i = I$, so $B_i = A_i^{-1}$.
Since it is only necessary to invert matrices of size $2 \times 2$, for simplicity we will content ourselves with using a direct matrix inverse to compute the new derivatives, though other methods could be used to compute the desired derivative and integral values.

Now, given the matrices $B_i$, we must now use the coefficients for the desired linear combinations stored as columns of $B_i$ to compute the corresponding integral terms.
Using similar reasoning as before, taking the needed linear combinations of the integral terms stored in \textit{ints} corresponds to right-multiplication of each $2\times 2$ matrix corresponding to a given degree and interval by the matrix $B_i$ corresponding to that interval.
This routine must return both the new scaled version of \textit{ints} and the corresponding matrices $B_i$ (these are the $i$'th entry along the first axis of \textit{ints}) because the matrices $B_i$ must be used again later to write the computed coefficients for the general functions in terms of the original $u$ and $v$ rather than their scaled linear combinations.
The internal workings of this auxiliary routine are outlined in Algorithm \ref{ScaleKnotFuncs}.

\begin{algorithm}
\caption{Take linear combinations of the input knot functions such that the desired linear combinations will satisfy the constraints $u_i\pr{t_i} = v_i\pr{t_{i+1}} = 1$ and $u_i\pr{t_{i+1}} = v_i\pr{t_i} = 0$.
Return the corresponding integral terms of these linear combinations and the coefficients for the linear combinations over each interval.}
\begin{algorithmic}[1]
\Procedure{ScaleKnotFuncs}{$ints, Tlens, tol=1E-8$}
	\LineComment{Copy \textit{ints} so it can be modified in-place without modifying the input array.}
	\State $ints = \text{copy}\pr{ints}$
	\LineComment{Get a boolean array showing where the the lengths of the intervals are nonzero.}
	\State $pos = Tlens > tol$
	\LineComment{Compute the coefficients of the needed linear combinations.}
	\State $invs = $ array of $0$'s of shape $\shape{ints}\br{1:}$
	\State $invs\br{pos} = \func{inv}{ints\br{0, pos}}$
	\LineComment{Perform matrix multiplication of each set of coefficients for each interval and degree}
	\LineComment{by the corresponding scaling matrix for each interval.}
	\State $ints\br{:} = \func{MatMul}{ints, invs}$
	\State \Return $ints, invs$
\EndProcedure
\end{algorithmic}
\label{ScaleKnotFuncs}
\end{algorithm}

Once the main loop in Algorithm \ref{gb_basis_construction} is finished, the computed general function coefficients must be changed to represent each basis function in terms of the original knot functions rather than the chosen linear combinations of them.
This is equivalent to left-multiplying the set of coefficients for each interval by the matrix $B_i$ (as in the explanation for \textit{ScaleKnotFuncs}).
This process is shown in Algorithm \ref{ScaleGenFuncCoefs}.

\begin{algorithm}
\caption{Perform a change of basis on the general function coefficients so that the general function coefficients used to represent the basis correspond to the functions originally used to form the array \textit{ints} of integral values.}
\begin{algorithmic}[1]
\Procedure{ScaleGenFuncCoefs}{$invs, genfunc$}
	\State \Return $\func{MatMul}{invs, genfunc\br{\dots,None}}\br{\dots,0}$
\EndProcedure
\end{algorithmic}
\label{ScaleGenFuncCoefs}
\end{algorithm}

The auxiliary routine \textit{PolyInt} is dependent on the polynomial representation used.
The array \textit{Tvals} is used as an argument because the polynomial basis used could be defined over some given interval (as are the Bernstein, Chebyshev, and Legendre polynomials).
For the power basis polynomials, that argument is not needed.
As has already been mentioned, the interval lengths are all that is necessary since the polynomial and general function terms are all assumed to be shifted to be defined over an interval starting at $0$.

The auxiliary routine \textit{PolyVal} should be handled similarly as \textit{PolyInt}.
This routine is also dependent on the polynomial representation and is easily defined as using Horner's algorithm, the De Casteljau algorithm, Clenshaw's algorithm, or some other polynomial evaluation algorithm.

\textit{GenFuncInts} is a function to compute the general function integrals
\[\int_{t_i}^{t_{i+1}} \pr{a_{i,j}^{d-2} u_{i,j}^{\br{d-2}}\pr{s} + b_{i,j}^{d-2} v_{i,j}^{\br{d-2}}\pr{s}} ds\]
with the corresponding constant terms
\[- a_{i,j}^{d-2} u_j^{[d-2]}(t_j) - b_{i,j}^{d-2} v_j^{[d-2]}(t_j)\]
from the left endpoint of the integral.
It is dependent on the representation used for the knot functions (we've only introduced using the knot functions themselves thus far).
In the case that the knot functions themselves are used, Algorithm \ref{GenFuncInt} shows how this can be done.

\begin{algorithm}
\caption{Integrate the general function part of each basis function over each interval in the support of that basis function.}
\begin{algorithmic}[1]
\Procedure{GenFuncInt}{$genfunc, wints$}
	\State $consts = - \func{sum}{genfunc * wints\br{:,:,0}, axis=-1}$
	\State $vals = \func{sum}{genfunc * wints\br{:,:,1}, axis=-1} + consts$
	\State \Return $vals, consts$
\EndProcedure
\end{algorithmic}
\label{GenFuncInt}
\end{algorithm}

\textit{IntegrateSupports} is a function that evaluates the integrals of each basis function over its corresponding support.
It should return both the desired indefinite integrals and the constant terms (stored in variable \textit{consts}) that come from the left bounds of each integral.
This function is outlined in Algorithm \ref{IntegrateSupports}.

\begin{algorithm}
\caption{Compute the definite integrals of each basis function over its support.}
\begin{algorithmic}[1]
\Procedure{IntegrateSupports}{$Tvals, pints, genfunc, wints$}
	\LineComment{\textit{wints} and \textit{Tvals} line the integral terms and the interval lengths up}
	\LineComment{with their corresponding interval in each basis function.}
	\State $vals, consts = \func{GenFuncInt}{genfunc, wints}$
	\State $vals += \func{PolyVal}{pints, Tvals, Tvals}$
	\State \Return $\func{sum}{vals, axis=-1}, consts$
\EndProcedure
\end{algorithmic}
\label{IntegrateSupports}
\end{algorithm}

\textit{OffsetConstants} is another helper routine that interfaces with the polynomials and is dependent on how the polynomials are represented.
It adds the terms stored in \textit{consts} to the corresponding terms in \textit{pints}.
In the case of the power basis, Chebyshev basis, or Legendre basis, this can be done by adding each constant term to the term in the polynomial representation that represents constants.
In the case of the Bernstein polynomials, since all the coefficients sum to $1$, adding a constant is the same as adding a constant to each coefficient, so this operation can be performed by adding the constant term for each polynomial to all the coefficients for that polynomial.

\textit{OffsetDifferences} is an auxiliary routine that takes care of differencing between the integrated terms from the previous basis function to form the differences over each interval that are needed to form the new basis.
Once this function has been applied, the terms account for everything included in the recurrence in Definition \ref{gbspline_standard_recurrence} with the exception of the constant term for each interval that comes from evaluating each basis function on the right endpoint of the interval to the left of the current interval.
This function also does not account for adding the ones to handle basis functions that are identically $0$.
The pseudocode for this algorithm is outlined in Algorithm \ref{OffsetDifferences}.

\begin{algorithm}
\caption{Take the differences between the integral terms for the previous set of basis functions to start forming the new set of basis functions.}
\begin{algorithmic}[1]
\Procedure{OffsetDifferences}{$pints, genfunc$}
	\State $n = \func{shape}{pints}\br{0}$
	\State $nints = \func{shape}{pints}\br{1}$
	\LineComment{Allocate the arrays needed to store the coefficients for the new basis.}
	\State $npolys = $ new array of $0$'s of shape $\pr{n-1, nints+1, dmin}$
	\State $ngenfunc = $ new array of $0$'s of shape $\pr{n-1, nints+1, 2}$
	\LineComment{Take the needed differences between the corresponding terms.}
	\State $npolys\br{:,:-1} \pluseq pints\br{:-1}$
	\State $npolys\br{:,1:} \mineq pints\br{1:}$
	\State $ngenfunc\br{:,:-1} \pluseq genfunc\br{:-1}$
	\State $ngenfunc\br{:,1:} \mineq genfunc\br{1:}$
	\State \Return \textit{npolys}, \textit{ngenfunc}
\EndProcedure
\end{algorithmic}
\label{OffsetDifferences}
\end{algorithm}

\textit{AddOnes} is used to add the ones that come from the integral terms from Definition \ref{gbspline_definition} that correspond to basis functions that are identically $0$.
We have separated it as an auxiliary routine because it both depends on the polynomial basis used.
We add $1$only to the last interval of basis functions for which the first term of the recurrence from Definition \ref{gbspline_definition} corresponds to a basis function that is identically $0$.
This is because, when constructing the basis functions of the next highest degree, the integral term corresponding to a basis function with index $i$ appears only in the expressions for the basis functions at index $i-1$ and $i$.
Of those two basis functions, only the basis function at index $i$ takes nonzero values on an interval that lies to the right of the support of the basis function that is identically $0$.
Once understood, this operation is very simple to perform, as can be seen in Algorithm \ref{AddOnes}, which demonstrates this auxiliary routine for polynomials represented in the power basis.
Though this routine depends on the polynomial representation used, it is not necessary to pass \textit{Tvals} since a constant terms is the same for a polynomial represented over any interval.

\begin{algorithm}
\caption{Add ones where needed to account for the integral terms in Definition \ref{gbspline_definition} that correspond to basis functions that are identically $0$.}
\begin{algorithmic}[1]
\Procedure{AddOnes}{$polys, pos$}
	\LineComment{Where the integral of the basis function of previous degree at the same}
	\LineComment{index was $0$, add $1$ to the constant term of the last polynomial term.}
	\State $polys\br{\sim pos\br{:-1},-1,-1} \pluseq 1$
\EndProcedure
\end{algorithmic}
\label{AddOnes}
\end{algorithm}

\textit{ConnectBoundaries} is the last auxiliary routine needed to construct the new basis functions from the previous ones.
In the recurrence in Definition \ref{gbspline_standard_recurrence}, this function adds in the terms $\Nof{i}{p}{t_j}$.
This function effectively starts at the leftmost interval in the support of each basis function, computes the value of the basis function at the end of that interval, adds that constant term to the polynomial term of the basis function on the next interval, and continues until it has added the needed constant terms to every interval in the support of that basis function.
This is done in a vectorized manner in Algorithm \ref{ConnectBoundaries}.

\begin{algorithm}
\caption{For each basis funciton, add in the constant terms $\Nof{i}{p}{t_j}$ to each interval where they are needed.}
\begin{algorithmic}[1]
\Procedure{ConnectBoundaries}{$polys, genfunc, wints, Tvals$}
	\LineComment{For each basis function, evaluate all but the leftmost polynomial term at the end}
	\LineComment{of the interval where it is defined.}
	\State $vals = \func{PolyVal}{polys\br{:,:-1}, Tvals, Tvals\br{:,:-1}}$
	\LineComment{Add in the corresponding values of the general functions.}
	\State $vals \pluseq \func{sum}{genfunc\br{:,:-1} * wints\br{:-1,:,1}, axis=-1}$
	\LineComment{Add the needed constants to their corresponding polynomial terms.}
	\State $\func{OffsetConstants}{polys\br{:,1:}, \func{cumsum}{vals, axis=1}}$
\EndProcedure
\end{algorithmic}
\label{ConnectBoundaries}
\end{algorithm}

\section{Conclusion}
\label{sec:conclusion}
In this manuscript we have presented an algorithm for the evaluation of GB-spline curves via their piecewise representation which is more direct than the recursive integral process given in the original definition for GB-splines.  The new algorithm makes practical computation simpler and easier to implement.  Moreover, using piecewise local representations make it so that the cost of evaluating a given spline curve is bound primarily by the costs of finding the portion of the knot vector in which a given point lies and evaluating the functions spanned by the spline basis. The computational routines here can also be used to work toward developing more efficient methods for the evaluation of specific classes of GB-spline curves.
They provide working examples that can be used to further study possible ways to provide better evaluation routines or subdivision methods for specific classes of GB-spline curves.


The use of piecewise local representations for the evaluation of GB-spline curves motivates the use of these local representations for other operations. One such operation is GB-spline refinement which will be covered in a future paper by the authors.  Also the local bases on each interval used to construct each spline curve share some of the useful properties of the bases used in \cite{scott_extraction}.
There, the process of inserting knots to represent a given B-spline curve as a piecewise polynomial was used to develop a local element structure that can, in turn, be used in isogeometric analysis~\cite{iga_gb_splines,Hughes2005}.
The local representations introduced here can be used in a similar manner to provide element structures for Isogeometric Analysis.

\bibliographystyle{elsarticle-num}
\bibliography{thesis}

\begin{thebibliography}{10}
\expandafter\ifx\csname url\endcsname\relax
  \def\url#1{\texttt{#1}}\fi
\expandafter\ifx\csname urlprefix\endcsname\relax\def\urlprefix{URL }\fi
\expandafter\ifx\csname href\endcsname\relax
  \def\href#1#2{#2} \def\path#1{#1}\fi

\bibitem{shapepreserv}
E.~Mainar, J.~Peña, J.~Sánchez-Reyes, Shape preserving alternatives to the
  rational b{e}ézier model, Computer Aided Geometric Design 18~(1) (2001) 37 --
  60.
\newblock \href
  {http://dx.doi.org/http://dx.doi.org/10.1016/S0167-8396(01)00011-5}
  {\path{doi:http://dx.doi.org/10.1016/S0167-8396(01)00011-5}}.

\bibitem{electromag}
A.~Buffa, G.~Sangalli, R.~Vázquez, Isogeometric analysis in electromagnetics:
  B-splines approximation, Computer Methods in Applied Mechanics and
  Engineering 199~(17?20) (2010) 1143 -- 1152.
\newblock \href {http://dx.doi.org/http://dx.doi.org/10.1016/j.cma.2009.12.002}
  {\path{doi:http://dx.doi.org/10.1016/j.cma.2009.12.002}}.

\bibitem{Ksasov1999}
B.~I. Ksasov, P.~Sattayatham, G{B}-splines of arbitrary order, Journal of
  Computational and Applied Mathematics 104~(1) (1999) 63 -- 88.
\newblock \href {http://dx.doi.org/10.1016/S0377-0427(98)00265-9}
  {\path{doi:10.1016/S0377-0427(98)00265-9}}.

\bibitem{costantini2005}
P.~Costantini, T.~Lyche, C.~Manni, On a class of weak {T}chebycheff systems,
  Numerische Mathematik 101~(2) (2005) 333--354.
\newblock \href {http://dx.doi.org/10.1007/s00211-005-0613-6}
  {\path{doi:10.1007/s00211-005-0613-6}}.

\bibitem{ue_spline_original}
G.~Wang, M.~Fang, Unified and extended form of three types of splines, Journal
  of Computational and Applied Mathematics 216~(2) (2008) 498 -- 508.
\newblock \href {http://dx.doi.org/10.1016/j.cam.2007.05.031}
  {\path{doi:10.1016/j.cam.2007.05.031}}.

\bibitem{iga_gb_splines}
C.~Manni, F.~Pelosi, M.~L. Sampoli, Generalized {B}-splines as a tool in
  isogeometric analysis, Computer Methods in Applied Mechanics and Engineering
  200~(5–8) (2011) 867 -- 881.
\newblock \href {http://dx.doi.org/10.1016/j.cma.2010.10.010}
  {\path{doi:10.1016/j.cma.2010.10.010}}.

\bibitem{exponential_subdivision}
L.~Romani, From approximating subdivision schemes for exponential splines to
  high-performance interpolating algorithms, Journal of Computational and
  Applied Mathematics 224~(1) (2009) 383 -- 396.
\newblock \href {http://dx.doi.org/10.1016/j.cam.2008.05.013}
  {\path{doi:10.1016/j.cam.2008.05.013}}.

\bibitem{subdivision_book}
J.~Warren, H.~Weimer, Subdivision Methods for Geometric Design: A Constructive
  Approach, Morgan Kaufmann series in computer graphics and geometric modeling,
  Morgan Kaufmann, 2002.

\bibitem{quasiinterpolation}
P.~Costantini, C.~Manni, F.~Pelosi, M.~L. Sampoli,
  \href{http://dx.doi.org/10.1016/j.cagd.2010.07.004}{Quasi-interpolation in
  isogeometric analysis based on generalized {B}-splines}, Comput. Aided Geom.
  Design 27~(8) (2010) 656--668.
\newblock \href {http://dx.doi.org/10.1016/j.cagd.2010.07.004}
  {\path{doi:10.1016/j.cagd.2010.07.004}}.
\newline\urlprefix\url{http://dx.doi.org/10.1016/j.cagd.2010.07.004}

\bibitem{tmeshes}
C.~Bracco, F.~Roman, Spaces of generalized splines over t-meshes, Journal of
  Computational and Applied Mathematics 294 (2016) 102 -- 123.
\newblock \href {http://dx.doi.org/http://dx.doi.org/10.1016/j.cam.2015.08.006}
  {\path{doi:http://dx.doi.org/10.1016/j.cam.2015.08.006}}.

\bibitem{scott_extraction}
M.~J. Borden, M.~A. Scott, J.~A. Evans, T.~J.~R. Hughes, Isogeometric finite
  element data structures based on {B}{\'e}zier extraction of {NURBS},
  International Journal for Numerical Methods in Engineering 87~(1-5) (2011)
  15--47.
\newblock \href {http://dx.doi.org/10.1002/nme.2968}
  {\path{doi:10.1002/nme.2968}}.

\bibitem{Hughes2005}
T.~Hughes, J.~Cottrell, Y.~Bazilevs, Isogeometric analysis: {CAD}, finite
  elements, {NURBS}, exact geometry and mesh refinement, Computer Methods in
  Applied Mechanics and Engineering 194~(39–41) (2005) 4135 -- 4195.
\newblock \href {http://dx.doi.org/10.1016/j.cma.2004.10.008}
  {\path{doi:10.1016/j.cma.2004.10.008}}.

\end{thebibliography}

\end{document}